\documentclass[11pt]{amsart}

\usepackage[margin=1in]{geometry}
\usepackage{amsthm,amssymb}
\usepackage{amsmath}
\usepackage{tikz}
\usepackage{enumerate}
\usepackage[colorlinks=true,linkcolor=blue,anchorcolor=blue,citecolor=blue,filecolor=blue,menucolor=blue,runcolor=blue,urlcolor=blue]{hyperref}
\usepackage{tikz-cd}
\usepackage{mathtools}
\mathtoolsset{showonlyrefs} % Only show equation numbers for referenced equations
\usepackage{mathrsfs}
\usepackage[noadjust]{cite}
\usepackage{placeins}
\usepackage{enumerate}

\newtheorem*{theorem*}{{\textbf{Theorem}}}

\newtheorem*{lemma*}{{\textbf{Lemma}}}

\newtheorem*{corollary*}{{\textbf{Corollary}}}

\theoremstyle{definition}

\theoremstyle{remark}

\newcommand{\ZZ}{\mathbb Z}
\newcommand{\QQ}{\mathbb Q}
\newcommand{\RR}{\mathbb R}
\newcommand{\CC}{\mathbb C}
\renewcommand{\AA}{\mathbb A} % \AA is taken by something.

\newcommand{\GG}{\mathbb G}

\newcommand{\G}{\mathcal G}

\renewcommand{\pmod}[1]{\ (\mathrm{mod}\ #1)}

\newcommand{\bparen}[1]{\left(#1\right)}

\begin{document}

\title{On trace-one generators of abelian cubic fields}

\author{Andrew O'Desky}
\address{
\parbox{0.5\linewidth}{
    Department of Mathematics\\
    Princeton University\\
    %Princeton, NJ USA
    }
}
\email{andy.odesky@gmail.com}

\date{July 14, 2024}

\begin{abstract}
    Let $K$ be a tamely ramified abelian cubic number field 
    with discriminant $D_K$. 
    We prove that the number of 
    trace-one monic integral polynomials with root field $K$ 
    and height $H$ is 
    equal to the number of ideals 
    in the quadratic field $\QQ(\sqrt{-3})$ 
    with norm $H^2 D_K^{-1/2}$. 
\end{abstract}

\maketitle

%\section{Introduction} %{*

Our recent work \cite{cubic} found a new connection 
    between cubic abelian number fields 
    and the arithmetic of the quadratic field 
    $\QQ(\sqrt{-3})$. 
This connection exists because 
    the multiplicative group of $\QQ(\sqrt{-3})$ 
    is a factor of the unit group in 
    the group algebra of the finite group $\ZZ/3\ZZ$. 
This note describes 
    another facet of this connection. 

Fix an abelian cubic number field $K$ 
    which is {tamely ramified over $\QQ$} 
    and let $F_K$ denote the set of polynomials of the form 
    $t^3 -t^2 + at + b \in \ZZ[t]$ 
    whose root field is isomorphic to $K$. 
Let $D_K$ denote the discriminant of $K$. 
For such polynomials we necessarily have $a \leq 0$ and 
    we set $H(t^3 -t^2 + at + b) = \sqrt{1-3a}$ 
    (``toric height''). 
The main result of this note is the following formula: 

\begin{theorem*} %{*
%Let $\lambda_K = \min \{H(f)^2 : f \in F_K\}$. 
%Consider the Dirichlet series 
%\begin{equation} %{*
%    Z_K(s) = \sum_{f\in F_K} H(f)^{-2s}.
%\end{equation} %*}
\begin{equation} %{*
    %Z_K(s) = 
%    \sum_{\substack{f=t^3-t^2+at+b \in \ZZ[t]\\ G_f\cong \ZZ/3\ZZ,\,K_f\cong K}}(1-3a)^{-s}=
    \sqrt{D_K}^{s}
\sum_{f\in F_K} H(f)^{-2s}=
    (1-3^{-s})\zeta_{\QQ(\sqrt{-3})}(s)
\end{equation} %*}
where $\zeta_{\QQ(\sqrt{-3})}$ is the Dedekind zeta function 
    of $\QQ(\sqrt{-3})$. 
\end{theorem*} %*}

This shows that 
    the number of polynomials in $F_K$ 
    of a given height is 
    \emph{essentially independent of $K$ and determined 
    by the arithmetic of $\QQ(\sqrt{-3})$}. 
This is illustrated in the table below 
    for the fields 
    $K_{49} = \QQ(\zeta_7)^+$ and $K_{169}=\QQ[t]/(t^3-t^2-4t-1)$. 

\begin{figure}[h!] %{*
\scriptsize
\renewcommand{\arraystretch}{1.3}
    \begin{minipage}{0.5\textwidth}
\begin{equation} %{*
\begin{array}{r|l}
H(f)^2 & f: K_f = K_{49}\\
    \hline
 7 \times 1& t^3 - t^2 - 2t + 1\\
 7 \times 4& t^3 - t^2 - 9t + 1\\
 7 \times 7& t^3 - t^2 - 16t + 29,\,\, t^3 - t^2 - 16t - 13\\
 7 \times 13& t^3 - t^2 - 30t + 43,\,\, t^3 - t^2 - 30t - 41\\
 7 \times 16& t^3 - t^2 - 37t + 29 \\
 7 \times 19& t^3 - t^2 - 44t + 127,\,\, t^3 - t^2 - 44t - 83\\
 7 \times 25& t^3 - t^2 - 58t - 13\\
 7 \times 28& t^3 - t^2 - 65t + 169,\,\, t^3 - t^2 - 65t - 167\\
 7 \times 31& t^3 - t^2 - 72t + 169,\,\, t^3 - t^2 - 72t - 41\\
 7 \times 37& t^3 - t^2 - 86t + 337,\,\, t^3 - t^2 - 86t - 251\\
 7 \times 43& t^3 - t^2 - 100t + 113,\,\, t^3 - t^2 - 100t - 181\\
\end{array}
\end{equation} %*}
    \end{minipage}%
    \begin{minipage}{0.5\textwidth}
\begin{equation} %{*
\begin{array}{r|l}
H(f)^2 & f: K_f = K_{169}\\
    \hline
 13 \times 1& t^3 - t^2 - 4t - 1\\
 13 \times 4& t^3 - t^2 - 17t + 25\\
 13 \times 7& t^3 - t^2 - 30t + 25,\,\, t^3 - t^2 - 30t - 53\\
 13 \times 13& t^3 - t^2 - 56t + 181,\,\, t^3 - t^2 - 56t + 25\\
 13 \times 16& t^3 - t^2 - 69t - 131\\
 13 \times 19& t^3 - t^2 - 82t + 155,\,\, t^3 - t^2 - 82t - 235\\
 13 \times 25& t^3 - t^2 - 108t + 337\\
 13 \times 28& t^3 - t^2 - 121t + 545,\,\, t^3 - t^2 - 121t - 79\\
 13 \times 31& t^3 - t^2 - 134t - 131,\,\, t^3 - t^2 - 134t - 521\\
 13 \times 37& t^3 - t^2 - 160t + 467,\,\, t^3 - t^2 - 160t - 625\\
 13 \times 43& t^3 - t^2 - 186t + 961,\,\, t^3 - t^2 - 186t + 415 
\end{array}
\end{equation} %*}
    \end{minipage}
\end{figure} %*}

\begin{corollary*}
%Let $O$ denote the ring of integers in $\QQ(\sqrt{-3})$. 
If $t^3 -t^2 + at + b \in F_K$ then $a \leq 0$ and 
    $\sqrt{D_K}$ divides $1-3a$. 
Fix $a \in \ZZ_{\leq 0}$. 
The number of polynomials of the form 
    $t^3 -t^2 + at + b \in F_K$ 
    for any $b \in \ZZ$ is equal to 
    the number of ideals in $\QQ(\sqrt{-3})$ 
    with norm $N=(1-3a)D_K^{-1/2}$. 
\end{corollary*}

\begin{figure} %{*
\scriptsize
\begin{equation} %{*
\begin{array}{|c|c|c|c|c|c|c|c|c|c|c|c|c|c|c|c|c|c|c|c|c|c|}
\hline
N&1&4&7&13&16&19&25&28&31&37&43&49&52&61&64&67&73&76&79&91&97 \\
\hline
d_N&1&1&2&2&1&2&1&2&2&2&2&3&2&2&1&2&2&2&2&4&2 \\
\hline
\end{array}
\end{equation} %*}
    \caption{Coefficients of 
    $\zeta_{\QQ(\sqrt{-3})}(s)=\sum_N d_N N^{-s}$ 
    for $N \equiv 1 \pmod 3$ up to $97$.}
\end{figure} %*}

From the corollary one easily derives the following formula 
for any $a \in \ZZ_{\leq 0}$: 
\begin{equation} %{*
    \#\{\,f=t^3 -t^2 + at + b \,\,:\,\, b \in \ZZ,\,\, K_f \cong K\,\}
    =
    \sigma_0\left(P_1\left(\frac{1-3a}{\sqrt{D_K}}\right)\right)
    %=\prod_{\substack{p |N\\p \equiv 1 \pmod 3}}
    %    (v_p(N)+1).
\end{equation} %*}
    where $\sigma_0(P)$ is the number of divisors of $P$ 
    and $P_1(N)$ is the largest divisor of $N$ 
    only divisible by primes $\equiv 1 \pmod 3$. 

\FloatBarrier

%Combining our theorem with the result from \cite{cubic} 
%    proves the following corollary. 
%
%\begin{corollary}
%\end{corollary}

\subsection*{Acknowledgements} %{*

%A.O. is grateful to 
This work was supported by NSF grant DMS-2103361. 
%The author is grateful to Manjul Bhargava for helpful discussions. 

%*}

%*}

\section{Proof of the theorem} %{*

\subsection{Poisson summation} %{*

%We use a refinement of the method from \cite{cubic}. 
Consider the rank $2$ algebraic torus 
    $\G = R^{\QQ(\sqrt{-3})}_\QQ \GG_m$ 
    and let $\AA$ denote the adele ring over $\QQ$. 
Let $f \in L^1(\mathcal G(\AA))$ and let $\widehat{f}$ 
    denote its Fourier transform. 
Let $\G(\QQ)^\perp$ denote the subgroup of characters in 
    the Pontryagin dual group 
    $\G(\AA)^\vee = \{\chi \colon \G(\AA) \to S^1\}$ 
    which are trivial on $\G(\QQ)$. 
The general Poisson summation formula --- 
following from the classical proof for $\ZZ \subset \RR$ --- 
says that if 
$\widehat{f}\,|_{\G(\mathbb Q)^\perp} \in L^1(\G(\mathbb Q)^\perp)$ 
then 
\begin{equation}\label{eqn:poissonae} %{*
    \int_{\G(\QQ)} f(xy)\, dx = 
    \int_{\G(\QQ)^\perp}
    \widehat{f}(\chi)
    \chi(y)
    \,d\chi 
\end{equation} %*}
for almost every $y \in \G(\mathbb A)$ and 
a suitably normalized Haar measure $d\chi$ on $\G(\QQ)^\perp$ 
\cite[Theorem~4.4.2, p.~105]{MR1397028}. 
The torus $\G$ may be identified with the group of units 
    of the group algebra of $G=\ZZ/3\ZZ$, 
    and therefore contains $G$ as a subgroup. 
We set $T = \G/G$, another rank $2$ torus. 
In \cite{cubic} this formula was applied with $y = 1$ 
and the torus $T$ in place of $\G$ 
to reexpress the height zeta function 
of a certain toric compactification of $T$. 

%*}

\subsection{Toric height} %{*

Batyrev--Tschinkel~{\cite{MR1369408}} introduced 
    a canonically defined height function 
    associated to any $T$-linearized line bundle 
    on a toric variety $S$ with open torus $T$. 
We recall some properties of this height 
    and set some notation. 
The group $\mathrm{Pic}_T(S)$ of $T$-linearized line bundles is 
    isomorphic to the free abelian group $\ZZ^{\Sigma(1)}$ on 
    the rays of the fan $\Sigma$ associated to $S$, 
    and we denote the toric height of $x \in T(\AA)$ 
    with respect to the toric divisor 
    $s = \sum_{e \in \Sigma(1)} s_eD_e$ 
    by $H(x,s)$. 
For fixed $s \in \ZZ^{\Sigma(1)}$, 
    the toric height extends to 
    the set of adelic points $T(\AA)$ 
    and decomposes as a product of local heights. 
Moreover, since the toric height is valued in 
    the positive reals, 
    it extends linearly to a pairing 
\begin{equation} %{*
    H \colon T(\AA) \times \CC^{\Sigma(1)} \to \CC.
\end{equation} %*}
For more details see \cite[\S3]{cubic}. 

%*}

\subsection{An orbit parametrization for \texorpdfstring{$\ZZ/3\ZZ$}{Z/3Z}-algebras} %{*

In \cite[\S2]{odesky2023moduli} it was shown that 
    there is a canonical bijection 
    between rational points of $T$ and 
    equivalence classes of pairs $(K,x)$ 
    where $x$ is a trace-one normal element of $K$. 
In \cite{cubic} the toric height of 
    any rational point $(K,x) \in T(\QQ)$ 
    whose characteristic polynomial is in $F_K$ 
    was computed explicitly: 
there is a toric divisor $D \subset S$ 
    defined over $\ZZ$ 
    such that any rational point $(K,x)$ 
    whose characteristic polynomial 
    $f_x=t^3-t^3+at+b$ is in $F_K$ is $D$-integral 
    \cite[Prop.~3]{cubic}, and  
\begin{equation} %{*
    H(t^3-t^3+at+b,D) = \sqrt{1-3a}. 
\end{equation} %*}

%\begin{remark}
%For each $f \in F_K$ there are precisely two rational points 
%    of $T(\QQ)$ with characteristic polynomial $f$ 
%    \cite[Example~1]{cubic}. 
%\end{remark}

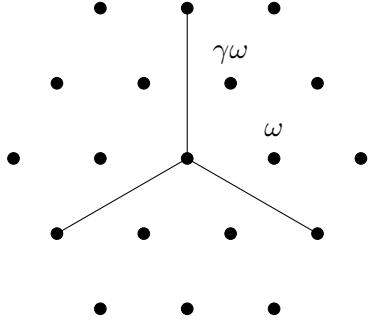
\begin{figure}[h!] %{*
    \centering
\begin{tikzpicture}[scale=2] 

\filldraw [black] (0,0) circle (.1em);
\filldraw [black] ({1/sqrt(3)},0) circle (.1em);
\filldraw [black] ({2/sqrt(3)},0) circle (.1em);
\filldraw [black] ({-1/sqrt(3)},0) circle (.1em);
\filldraw [black] ({-2/sqrt(3)},0) circle (.1em);
\filldraw [black] ({1/sqrt(3)+1/(2*sqrt(3))},1/2) circle (.1em);

\filldraw [black] ({1/(2*sqrt(3))},1/2) circle (.1em);
\filldraw [black] ({2/(2*sqrt(3))},1) circle (.1em);
\filldraw [black] ({-1/(2*sqrt(3))},-1/2) circle (.1em);
\filldraw [black] ({-2/(2*sqrt(3))},-1) circle (.1em);
\filldraw [black] ({-1/(2*sqrt(3))},1/2) circle (.1em);
\filldraw [black] ({-2/(2*sqrt(3))},1) circle (.1em);
\filldraw [black] ({1/(2*sqrt(3))},-1/2) circle (.1em);
\filldraw [black] ({2/(2*sqrt(3))},-1) circle (.1em);

\filldraw [black] ({sqrt(3)/2},-1/2) circle (.1em); % e1
\filldraw [black] ({-sqrt(3)/2},-1/2) circle (.1em);
\filldraw [black] (0,1) circle (.1em);

\filldraw [black] ({-sqrt(3)/2},1/2) circle (.1em);
\filldraw [black] (0,-1) circle (.1em);
 
\draw (0,0) -- ({sqrt(3)/2},-1/2);
\draw (0,0) -- (0,1);
\draw (0,0) -- ({-sqrt(3)/2},-1/2);

%\draw[dashed] (0,0) -- ({2/sqrt(3)},0); % w1
%\draw[dashed] (0,0) -- ({2/(2*sqrt(3))},1); % w2
%\draw[dashed] (0,0) -- ({-2/(2*sqrt(3))},1); % w3
%\draw[dashed] (0,0) -- ({-2/sqrt(3)},0); % w4
%\draw[dashed] (0,0) -- ({-2/(2*sqrt(3))},-1); % w5
%\draw[dashed] (0,0) -- ({2/(2*sqrt(3))},-1); % w6

\node[yshift=1em] at ({1/sqrt(3)},0) {$\omega$};
\node[yshift=1em] at ({1/(2*sqrt(3))},1/2) {$\gamma \omega$};
%\node[yshift=1em] at ({sqrt(3)/2},-1/2) {$v_1$};
%\node[yshift=1em] at (0,1) {$v_2$};
%\node[yshift=1em] at (-{sqrt(3)/2},-1/2) {$v_0$};
%\node[yshift=-1em] at (0,0) {$0$};

%\foreach \x in {-1,...,1}
%\foreach \y in {-1/2,1/2} 
%{
%   %\draw (\x,\y) circle (1.5pt); 
%   \filldraw [black] (\x,\y) circle (1.5pt); 
%}

\end{tikzpicture}  %*}
\vspace{.5cm}
\caption{The fan $\Sigma$ of the toric surface $S$ 
    with Galois conjugation $\omega \mapsto \gamma \omega$.}
\end{figure}%*}

Let $x$ be any trace-one normal element of $K$, 
    and let $(K,x)$ denote the associated rational point 
    of $T(\QQ)$ \cite[Theorem~4]{cubic}. 
We will apply the Poisson formula 
for a certain fixed $y \in \G(\AA)$ to the function 
\begin{equation} %{*
    f(v) = H(v(K,x),-s) 1_{D}
\end{equation} %*}
where $1_{D} \colon T(\AA) \to \{0,1\}$ 
    is the indicator function for $D$-integral points. 

    %*}

\subsection{Normal integral bases} %{*

The only dependence on $K$ in the Poisson formula 
    will be in our choice of 
    the ``twisting parameter'' $y \in \G(\AA)$. 
This parameter will be obtained as an approximation 
to the rational point $(K,x) \in T(\QQ)$. 
To construct it, we use a classical theorem of Noether 
    on the existence of normal integral bases. 

\begin{lemma*}
%Assume $K/\QQ$ is tamely ramified. 
There is an element $y=y_K \in \G(\AA)$ 
    and an element $k$ of the maximal compact subgroup 
    of $T(\AA)$ 
    such that $(K,x) = yk$.
    In particular, $H((K,x),s) = H(\pi(y),s)$ 
    where $\pi \colon \G \to T$ denotes 
    the natural quotient morphism. 
\end{lemma*}

\begin{proof} %{*
Let $v$ be a finite place of $\QQ$. 
A classical theorem of Noether says that $O_v/\ZZ_v$ 
has a normal integral basis if and only if 
    $K_v/\QQ_v$ is tamely ramified. 
This implies that 
\begin{equation} %{*
    (K_v/\QQ_v,x) = y_v k_v
\end{equation} %*}
    for some element $k_v$ of the maximal compact subgroup 
    of $T(\QQ_v)$. 
Since $K/\QQ$ is split at infinity, 
    the real point $(K_\RR/\RR,x)$ is equal to $\pi(y_\infty)$ 
    for some $y_\infty \in \G(\RR)$. 
Now take $y = (y_v)_v$ and $k =(k_v)$. 
\end{proof} %*}

%*}

\subsection{Proof of the theorem} %{*

Write $H(-,s,D) \coloneqq H(-,s)1_{D}$.  
By the lemma there is an element $y_K \in \G(\AA)$ such that 
    $$H(v(K,x),-s,D) = H(v\pi(y_K),s,D)$$ 
    for any $v \in \G(\QQ)$. 
The Poisson formula 
    for $f(v) = H(v(K,x),-s,D)$ 
    with $y = y_K$ implies that 
\begin{equation}\label{eqn:twistedPoisson} %{*
    \sum_{v \in \G(\QQ)}
        H(v(K,x),-s,D)
    =
    \int_{\G(\QQ)^\perp}
    \widehat{f}(\chi)\chi(y_K)\,d\chi.
\end{equation} %*}
Let $Z_K$ denote the left-hand side of the equality 
    in the theorem. 
There are precisely two rational points on $T$ corresponding 
    to a given polynomial with Galois group $\ZZ/3\ZZ$, 
    namely $(K,x)$ and $(K',x)$ where $K'$ is the twist of 
    the $\ZZ/3\ZZ$-algebra $K$ 
    by the outer automorphism of $\ZZ/3\ZZ$ 
    \cite[Example~1]{cubic}, 
    so the left-hand side of \eqref{eqn:twistedPoisson} 
    is equal to $2 Z_K$. 

To evaluate the right-hand side of \eqref{eqn:twistedPoisson} 
    we make use of the calculation of 
    $\widehat{f}(\chi)$ from \cite[Prop.~4]{cubic} 
    maintaining the same notation. 
In addition to the new factor $\chi(y_K)$ in the integrand, 
    there are a few differences since we are using $\G$ 
    rather than $T$ to compute the Fourier transform. 
At each place, we must compute the Fourier transform 
    with respect to the split locus --- 
    the image of $\G(\QQ_v) \to T(\QQ_v)$ --- 
    rather than $T(\QQ_v)$. 
At finite places $v \neq 3$ this results in 
the same expression for $\widehat{H_v}$ 
    except we instead sum over 
    the simpler sublattice (writing $E = \QQ(\sqrt{-3})$) 
\begin{equation} %{*
    \ZZ\langle v_1, v_2 \rangle 
    = X_\ast(\G_E)^{\Gamma(w/v)}
    \subset 
    X_\ast(T_E)^{\Gamma(w/v)} = \ZZ \langle v_1,\omega \rangle
\end{equation} %*}
where $w$ is any place of $E$ over $v$, 
$\omega = \tfrac23 v_1 + \tfrac13 v_2$ 
and $v_1,v_2$ are the cocharacters of $\G_E$ 
corresponding to the two nontrivial representations of $\ZZ/3\ZZ$, 
and $\Gamma(w/v)$ is the decomposition group of 
the Galois group of $E/\QQ$ at $w$. 
%This sublattice is the image of the split locus at $v$ modulo 
%    the maximal compact subgroup. 
The Fourier transform at $v = \infty$ is unchanged 
since 
    $X_\ast(\G_E)_\RR^{\Gamma(w/\infty)}= 
    X_\ast(T_E)_\RR^{\Gamma(w/\infty)}$. 
At $v = 3$ the local height is compactly supported, 
    and the computation in \cite{cubic} 
    shows that it does not contribute 
    to the right-hand side of \eqref{eqn:twistedPoisson}. 
This results in the following integral representation for $Z_K$: 
\begin{equation} %{*
    \int_{\G(\QQ)^\perp}
    \widehat{f}(\chi)\chi(y_K)\,d\chi
    =
    \left(\frac{-1}{2\pi i}\right)
    \frac{s}{\pi i}
\sum_{\eta}
    \eta^{-s}
\int_{\RR} 
    \frac{\chi_t(\eta)^{-1}\chi_t(y_K)\,dt}
    {(t+\frac{s}{\pi i})(t-\frac{s}{2\pi i})}
\end{equation} %*}
where now we only sum over \emph{split} $\eta$ 
corresponding to the image of $\G(\AA^f) \to T(\AA^f)$ 
on finite adeles; 
explicitly, $\eta = (m_w)_{v \neq 3,\infty}$ 
where 
\begin{equation} %{*
    m_w \in 
    \begin{cases}
        \ZZ_{\geq 0} \langle v_1,v_2 \rangle& \text{if $w$ split,}\\
        \ZZ_{\leq 0} \langle v_0\rangle & \text{otherwise.}\\
    \end{cases}
\end{equation} %*}
The integral can be evaluated by Cauchy's residue formula, 
and summing over $\eta$ obtains 
\begin{equation} %{*
    Z_K = 
    \chi_{\frac{s}{\pi i}}(y_K)^{-1}
    \bparen{\prod_{q\equiv2\pmod{3}}
    \sum_{c_q= 0}^\infty 
        q^{-2c_qs}}
    \bparen{
        \prod_{p\equiv1\pmod{3}}
        \sum_{a_p,b_p= 0}^\infty
            p^{-(a_p+b_p)s}} . 
\end{equation} %*}
This evaluates to 
$\chi_{\frac{s}{\pi i}}(y_K)^{-1}(1-3^{-s})
\zeta_{\QQ(\sqrt{-3})}$. 

In particular we see that 
$\chi_{\frac{s}{\pi i}}(y_K)^{-1}
=\lambda_K^{-s}$ where $\lambda_K = \min \{H(f)^2 : f \in F_K\}$. 
It remains to be seen that $\lambda_K = \sqrt{D_K}$; 
this follows from an elementary computation as follows. 
Let $|\!|\cdot|\!|$ denote the canonical norm on 
Minkowski space $K_\RR$. 
Let $|\!|\cdot|\!|_1$ denote 
the quotient norm on $K_\RR/\RR 1$ given by 
    $|\!|v+\RR1|\!|_1\coloneqq 
    \min \{ |\!|v-r1|\!| : r \in \RR\}$. 
Writing $v = (x,y,z) \in K_\RR$, one finds that 
%the minimum of $r\mapsto |\!|v-r1|\!|^2=p_2-2re_1+3r^2$ is achieved 
%when $r = \tfrac13$, so 
$|\!|v+\RR1|\!|_1^2 
    =\tfrac13(3p_2-1)=\tfrac23(1-3e_2)$, 
so the toric height and the quotient norm 
    are proportional by $\sqrt{\tfrac23}$. 
The group $\ZZ/3\ZZ$ acts unitarily on 
$K_\RR/\RR 1$ for $|\!|\cdot|\!|_1$, 
so the lattice $O_K/\ZZ1$ is the standard hexagonal lattice 
which has covolume $\frac{\sqrt3}{2} \ell^2$ 
where $\ell$ is the minimal length of a nonzero lattice point. 
Thus $\ell = \sqrt{\tfrac23} \sqrt{\lambda_K}$. 
On the other hand, the covolume of $O_K/\ZZ1$ 
is $\sqrt{3}^{-1} \cdot \mathrm{covol}(O_K \subset K_\RR) 
= \sqrt{3}^{-1} \sqrt{D_K}$, 
and equating both formulas for the covolume of $O_K/\ZZ1$ 
shows that $\sqrt{D_K} = \lambda_K$. 

%*}

%*}

\bibliography{draft}
\bibliographystyle{abbrv}

\end{document}